\documentclass[twoside,12pt]{article}
\usepackage[all]{xy}
\usepackage{amssymb,amsmath,latexsym}

\usepackage{amsfonts}
\usepackage{ulem}
\usepackage{amscd}
\usepackage{amsxtra}
\usepackage{mathrsfs}

\newtheorem{theorem}{Theorem}[section]
\newtheorem{lemma}[theorem]{Lemma}

\newtheorem{question}{\sc Question}

\newtheorem{corollary}[theorem]{Corollary}

\newtheorem{remark}{\sc Remark}

\numberwithin{equation}{section}

\newcommand{\cqfd}
{\hspace{1cm}
\rule{2mm}{2mm}%
\medbreak%
\par%
}

\def\pr{{\parindent0pt {\bf Proof.\ }}}
\def\cqfd
{\hspace{1cm}
\rule{2mm}{2mm}%
\medbreak%
\par%
}

\def\tri{{\rm Tri}}

\author{}

\begin{document}
\title{$f$-Biderivations and Jordan biderivations of unital algebras with idempotents}

\date{}
\maketitle \vspace*{-1.5cm}

\thispagestyle{empty}


\begin{center}
\author{Mohammad Ali Bahmani $^{1}$, Driss Bennis$^{2}$,  Hamid Reza Ebrahimi Vishki$^{3}$, and Brahim Fahid$^{4}$  }\end{center}

\noindent{\large\bf Abstract.}  The notion of $f$-derivations was introduced by Beidar and  Fong   to unify several  kinds of linear maps including derivations, Lie derivations and Jordan derivations. In this paper we introduce the notion of   $f$-biderivations as a natural ``biderivation" counterpart of the notion of ``$f$-derivations". We first show, under some conditions, that any $f$-biderivation is a Jordan biderivation. Then, we turn to study    $f$-biderivations of a unital algebra with an idempotent. Our second main result shows, under some conditions, that
 every Jordan  biderivation can be written as a sum of a biderivation, an antibiderivation and an extremal biderivation.   As a consequence we show that  every Jordan biderivation on a triangular algebra is a biderivation.   \bigskip




\section{Introduction}
Throughout the paper $\mathcal{R}$ will denote a commutative ring with unity,  $\mathcal{A}$ will be a unital $\mathcal{R}$-algebra with center $Z(\mathcal{A})$   and $\mathcal{M}$ will be a unital $\mathcal{A}$-bimodule.\medskip

Recall that an $\mathcal{R}$-linear map  $\mathcal{D}$ from $\mathcal{A}$ into $\mathcal{M}$ is said to be a derivation  (resp., an  antiderivation) if, for all $a,b \in \mathcal{A}$,  $  \mathcal{D}(ab)= \mathcal{D}(a)b+a\mathcal{D}(b)$ (resp., $ \mathcal{D}(ab)= \mathcal{D}(b) a+b\mathcal{D}(a)$). The inner derivations are classical examples of derivations. Recall that an $\mathcal{R}$-linear map  $\mathcal{D}$ is said to be inner if it is of the form  $\mathcal{D}(a)=[m_{0},a] $  for some $m_{0} \in \mathcal{M}$, where $[-,- ]$  stands for the Lie bracket.\\

In \cite{BF}, Beidar and  Fong introduced   the notion of $f$-derivations which unifies several particular kinds of linear maps  including the classical derivations as follows: \\
Consider a fixed nonzero multilinear polynomial $f$ in noncommuting indeterminates $x_{i}$ over $\mathcal{R}$:
\begin{equation}\label{eq-1}
  f(x_{1},\ldots,x_{n})=\sum_{\pi\in S_{n}}\alpha_{\pi}x_{\pi(1)}x_{\pi(2)}\ldots x_{\pi(n)},\quad \alpha_{\pi}\in \mathcal{R},
\end{equation}
 where $S_{n}$ denotes the symmetric group of order an integer $ n \geq 2 $. An $\mathcal{R}$-linear map $\mathcal{D} : \mathcal{A}\longrightarrow \mathcal{M} $ is called  an $f$-derivation  if  it satisfies
\begin{equation}\label{eq-2}
  \mathcal{D}(f(x_{1},\ldots,x_{n}))=\sum_{i=1}^{n}f(x_{1},\ldots,x_{i-1},\mathcal{D}(x_{i}),x_{i+1},\ldots, x_{n})
\end{equation}
for all $ x_{1},\ldots,x_{n} \in \mathcal{A}$.\\
Thus,
\begin{itemize}
    \item  a derivation  is an $f$-derivation for the polynomial $f(x_{1},x_{2}) =x_{1}x_{2}$,
   \item  a Jordan derivation is an $f$-derivation for the polynomial $ f (x_{1}, x_{2}) =x_{1}\circ x_{2}  := x_{1}x_{2} + x_{2}x_{1}$,
   \item a Jordan triple derivation is an $f$-derivation for the polynomial $ f(x_{1},x_{2},x_{3}) =x_{1}x_{2}x_{3} + x_{3}x_{2}x_{1}$,
   \item a Lie derivation is an $f$-derivation for the polynomial $ f (x_{1}, x_{2}) =[x_{1}, x_{2}]: = x_{1}x_{2} - x_{2}x_{1}$, and
   \item a Lie triple derivation is an $f$-derivation for the polynomial $ f(x_{1},x_{2},x_{3}) = [[x_{1},x_{2}] ,x_{3}]$.
\end{itemize}

In  \cite[Theorem 1.3]{B2},  Benkovi\v{c}   proved (under some conditions) that every $f$-derivation is a Jordan derivation. This means that in some situations studying  $f$-derivation is based on the study of   Jordan derivations. In \cite{BS},  Benkovi\v{c} and   \v{S}irovnik investigated  Jordan derivations on algebras with an idempotent. They proved that under certain ``nice" conditions every Jordan derivation  is a sum of a derivation and an antiderivation.\\

Our aim in this paper is to investigate the ``biderivation" counterpart of the above results. \\

Naturally one can define   a ``biderivation" counterpart of the $f$-derivations as follows: \\

In what follows, we consider  a fixed nonzero multilinear polynomial $f$ as defined in (\ref{eq-1}). An $\mathcal{R}$-linear map $ F: \mathcal{A}\times \mathcal{A} \longrightarrow \mathcal{M} $ is called an $f$-biderivation, if
 $$ F(f(x_{1},\ldots,x_{n}),z)=  \\
  \sum_{i=1}^{n}f(x_{1},\ldots,x_{i-1},F(x_{i},z),x_{i+1},\ldots, x_{n})  $$ and $$ F(z,f(x_{1},\ldots,x_{n}))=  \\
  \sum_{i=1}^{n}f(x_{1},\ldots,x_{i-1},F(z,x_{i}),x_{i+1},\ldots, x_{n})  $$
for all $ x_{1},\ldots,x_{n},z \in \mathcal{A}$.\\

Then,
\begin{itemize}
    \item  every $f$-biderivation $F$ is a biderivation when $f(x,y)=xy$   (see \cite{BrMM}),
   \item  every $f$-biderivation $F$ is a Jordan biderivation when  $f(x,y)=x\circ y$ (see for instance   \cite{AL}), and
   \item every $f$-biderivation $F$ is a Jordan triple biderivation when  $f(x,y,z)=xyz+zyx$ (see for example \cite{B89}).
\end{itemize}

We start our paper with the first main result, Theorem \ref{thm-lem}, which is the  ``biderivation" counterpart of Benkovi\v{c}'s result \cite[Theorem 1.3]{B2}.  It shows, under some conditions, that any $f$-biderivation is a Jordan biderivation. Then, in the remainder of the paper we focus  on the study of Jordan biderivations. Namely, we aim to establish the ``biderivation" counterpart of Benkovi\v{c} and    \v{S}irovnik's main result  \cite[Theorem 4.1]{BS}. As a main result (Theorem \ref{theorem-princ1}), we show, under some conditions, that  every Jordan  biderivation can be written as a sum of a biderivation, an antibiderivation and an extremal biderivation. Recall that a bilinear map $ D : \mathcal{A}\times \mathcal{A} \longrightarrow \mathcal{A}$ is called an antibiderivation if it is an antiderivation  with respect to both components.
A bilinear map $ D : \mathcal{A}\times \mathcal{A} \longrightarrow \mathcal{A}$ is called   an extremal biderivation if it is of the form $D(x, y) = [x, [y, a]]$ for all $x, y \in \mathcal{A}$, where $a \notin Z(\mathcal{A})$ and $[[\mathcal{A},\mathcal{A}], a] = 0$.\\

As a consequence of our second main result, we show that every Jordan biderivation of a triangular algebra is a biderivation (Corollary \ref{theorem-princ2}). Recall, for  two  $\mathcal{R}$-algebras  $A$ and  $B$  and an $(A,B)$-bimodule $M$, the  set   $$ \tri(A;M;B):=\{\left(                                                                                                                                      \begin{array}{cc}                                                                                                                                      a & m \\                                                                                                                                0 & b                                                                                                                               \end{array}                                                                                                                              \right) \mid a\in A, m\in M, b\in B\}$$
  equipped with the usual matrix operations is an $\mathcal{R}$-algebra called a (generalized) triangular $\mathcal{R}$-algebra  (see \cite{GJS}  for more details about this construction). In this paper    we assume that $M$ is also a faithful $ (A, B)$-bimodule. As interesting examples of   triangular matrix algebras one can cite the  (classical) upper triangular matrix algebras, the block upper triangular matrix algebras and the  nest algebras. It is important to recall  that  an algebra  $\mathcal{A}$ is isomorphic to a   triangular matrix algebra     if  there exists a non trivial idempotent $ e \in \mathcal{A}$ such that $ (1-e)\mathcal{A}e=0 $ (see, for instance, \cite[Theorem 5.1.4]{GJS}). Namely, in this case, $\mathcal{A}$ is isomorphic to $\tri(e\mathcal{A}e; e\mathcal{A}(1-e); (1-e)\mathcal{A}(1-e))$.

\section{Main results}

Let us start with the first main result  which  investigates $f$-biderivations  under some conditions. \medskip

 We say that an element $ r\in \mathcal{R}$ is $\mathcal{M}$-regular if  for every $m\in \mathcal{M}$,
$rm = 0$  implies $ m = 0$. Let $$\alpha=\sum_{\pi\in S_{n}}\alpha_{\pi} \in \mathcal{R}$$
be the sum of coefficients of the polynomial $f$ from (\ref{eq-1}).\\

The following result is  the ``biderivation" counterpart of \cite[Theorem 1.3]{B2}.

\begin{theorem} \label{thm-lem}
  Let $\mathcal{A}$ be a unital algebra and $\mathcal{M}$ a unital $\mathcal{A}$-bimodule. Let $ F : \mathcal{A}\times \mathcal{A} \longrightarrow \mathcal{M}$ be an $f$-biderivation,    with $ \alpha \neq 0$. If  $ \mathcal{M}$ is $(n-1)$-torsion free and $\alpha$ is  $\mathcal{M}$-regular, then $F$ is a  Jordan biderivation.
\end{theorem}
\pr First we prove that $F(1,y) = 0$ for all $y \in \mathcal{A}$. Let $ x_i = 1$ for $ i = 1, \ldots, n$. Then, by the definition of $f$-biderivation,
 \begin{equation*}
 \alpha F (1,y) = n\alpha F (1,y).
\end{equation*}
Then,  $(n -1) \alpha F (1,y) = 0$, and consequently $F (1,y) = 0$.

Now, we decompose the sum $f(x_1, x_2,   \ldots , x_n)=\sum_{\pi\in S_{n}}\alpha_{\pi}x_{\pi(1)}x_{\pi(2)}\ldots x_{\pi(n)}$ according to the order of $x_1$ and $x_2$ in the products $x_{\pi(1)}x_{\pi(2)}\ldots x_{\pi(n)}$. So let us decompose  $S_{n}$ into the following two disjoints subsets: $$S_{n}^{<}=\{\pi \in S_n; \pi^{-1}(1)< \pi^{-1}(2)\}\ \textrm{and} \ S_{n}^{>}=\{\pi \in S_n; \pi^{-1}(1) > \pi^{-1}(2)\}.$$
 Then, $f$ can be decomposed as a sum of $f^{<}$ and  $f^{>}$, where $$f^{<}(x_1, x_2,   \ldots , x_n)=\sum_{\pi\in S_{n}^{<}}\alpha_{\pi}x_{\pi(1)}x_{\pi(2)}\ldots x_{\pi(n)}$$
and $$f^{>}(x_1, x_2,   \ldots , x_n)=\sum_{\pi\in  S_{n}^{>}}\alpha_{\pi}x_{\pi(1)}x_{\pi(2)}\ldots x_{\pi(n)}.$$
It is clear that  $$f^{<}(x_1, x_2,  1, \ldots , 1)=\beta x_1x_2\ \textrm{where} \ \beta= \sum_{\substack{\pi\in \mathcal{S}_{n}^{<}}}\alpha_\pi $$
and  $$f^{>}(x_1, x_2,  1, \ldots , 1)= \gamma x_2x_1 \ \textrm{where} \  \gamma= \sum_{\substack{\pi\in \mathcal{S}_{n}^{>}}}\alpha_\pi.$$
Then, $f (x_1, x_2, 1, \ldots , 1)=\beta x_1x_2 + \gamma x_2x_1$. Since $F(1,y) = 0$,
\begin{equation*}
F(f (x_1, x_2, 1, \ldots , 1),y) = f (F(x_1,y) , x_2, 1, \ldots, 1) + f(x_1, F (x_2,y) , 1,\ldots , 1) .
\end{equation*}
Then,  for  all $x_1, x_2,y \in \mathcal{ A}$,
\begin{equation}\label{eq-4}
  F(\beta x_1x_2 + \gamma x_2x_1,y) = \beta F (x_1,y) x_2 + \beta x_1F (x_2,y) + \gamma F(x_2,y) x_1 + \gamma x_2F(x_1,y).
\end{equation}
Now we exchange the roles of $x_1$ and $x_2$ in (\ref{eq-4}) so that we get,  for  all $x_1, x_2,y \in \mathcal{ A}$,
\begin{equation}\label{eq-5}
  F(\beta x_2x_1 + \gamma x_1x_2,y) = \beta F (x_2,y) x_1 + \beta x_2F (x_1,y) + \gamma F(x_1,y) x_2 + \gamma x_1F(x_2,y).
\end{equation}
The sum of (\ref{eq-4}) and (\ref{eq-5}) is equal to
\begin{equation*}
F(\alpha x_1x_2 + \alpha x_2x_1,y) = \alpha F (x_1,y) x_2 + \alpha x_1F (x_2,y) + \alpha F (x_2,y) x_1 + \alpha x_2F (x_1,y)
\end{equation*}
for  all $x_1, x_2,y \in \mathcal{ A}$.
Since $\alpha$ is $ \mathcal{M}$-regular, we have,
for all $x_1, x_2,y \in \mathcal{A}$,
\begin{equation*}
F(x_1x_2 + x_2x_1,y) =  F (x_1,y) x_2 +  x_1F (x_2,y) +  F (x_2,y) x_1 +  x_2F (x_1,y).
\end{equation*}
 Similarly we prove that
\begin{equation*}
F(y,x_1x_2 + x_2x_1) =  F (y,x_1) x_2 +  x_1F (y,x_2) +  F (y,x_2) x_1 +  x_2F (y,x_1)
\end{equation*}
for all $x_1, x_2,y \in \mathcal{A}$. Therefore, $F$ is a  Jordan biderivation.
\cqfd

Now we turn to our second  aim of this paper. We study Jordan biderivations of unital algebras with idempotents.\\

Throughout the remainder of this section, we will fix the following condition  and notation:\\

\noindent \textbf{Setup  and  notation.} We assume that the algebra $\mathcal{A}$ admits  a nontrivial idempotent $e$. Then,
$$\mathcal{A}  = e\mathcal{A}e+e\mathcal{A}e'+e'\mathcal{A}e+e'\mathcal{A}e',$$  where $e'=1-e$. To simplify notation  we will use the following convention:\\
$a=eae\in e\mathcal{A}e=\mathcal{A}_{11}$, $m=eme'\in e\mathcal{A}e'=\mathcal{A}_{12}$, $n=e'ne\in e'\mathcal{A}e=\mathcal{A}_{21}$ and $b=e'be'\in e'\mathcal{A}e'=\mathcal{A}_{22}$.\\
Then each element $x=exe+exe'+e'xe+e'xe'\in \mathcal{A}$ can be represented in the form $x= exe+eme'+e'ne+e'be'=a+m+n+b$, where $a\in \mathcal{A}_{11}$, $m\in \mathcal{A}_{12}$, $n\in \mathcal{A}_{21}$ and $b\in \mathcal{A}_{22}$. Hence, every bilinear mapping
$J:\mathcal{A}\times \mathcal{A}\longrightarrow \mathcal{A}$  can be represented in the form
\begin{multline}
\label{linear3}
J(x,y) = J(a,a')+J(a,b')+J(a,m')+J(a,n')+J(m,a')+J(m,b')+J(m,m')\\
  +J(m,n')
+J(n,a')+J(n,b')+J(n,m')+J(n,n')+J(b,a')+J(b,b')\\
 +J(b,m')+J(b,n').
\end{multline}
for all $x=a+m+n+b, y=a'+m'+n'+b'\in  \mathcal{A}$.\\
Also,  in the rest of this paper we   assume  that any algebra, in particular  $\mathcal{A}$,  is $2$-torsion free (i.e.,  for every $x\in\mathcal{A}$,  $2x=0$ implies $x=0$).  Notice that in this case a bilinear map $ J : \mathcal{A}\times \mathcal{A} \longrightarrow \mathcal{A}$   is a Jordan biderivation if and only if, for all $x,y\in \mathcal{A}$, $J(x^2,y)=xJ(x,y)+J(x,y)x$  and $J(x,y^2)=yJ(x,y)+J(x,y)y$.\\


The second main result uses the following lemmas.

 \begin{lemma}\label{lemm-jbider}
    Let  $J : \mathcal{A} \times \mathcal{A }\longrightarrow \mathcal{A} $ be a Jordan biderivation. Then $[[x, y], J(x, y)] = 0$ for all $x, y \in \mathcal{A}$.
 \end{lemma}
\pr   Since $J$ is a Jordan  biderivation  we have for every $x,y\in \mathcal{A}$, $$ J(x^2,y^2)=xJ(x,y^2)+J(x,y^2)x.$$
  Then, $$J(x^2,y^2)=xyJ(x,y)+xJ(x,y)y+yJ(x,y)x+J(x,y)yx.$$
   By the same argument and using the fact that $J(x,y^2)=yJ(x,y)+J(x,y)y$   we get $$J(x^2,y^2)=yxJ(x,y)+yJ(x,y)x+xJ(x,y)y+J(x,y)xy.$$ Comparing both relations leads to  $[[x,y],J(x,y)]=0$.
\cqfd

\begin{lemma}\label{lemm-jbider2}
    Let  $ J : \mathcal{A} \times \mathcal{A} \longrightarrow \mathcal{A}$ be a Jordan biderivation. Then $J = J_1 + J_2$, where $J_1(x, y) = [x, [y, J(e, e)]]$ is an  extremal biderivation and $J_2$ is a Jordan biderivation such that $J_2(e, e) = 0$.
\end{lemma}
\pr
Let us first consider the identity $[[x,y],J(x,y)]=0$  for all $x,y\in \mathcal{A}$. Replacing  $x$ by $x+e$,  we get $$[[x,y],J(e,y)]+[[e,y],J(x,y)]=0.$$
 This implies that, $[[x,e],J(e,e)]=0$. Similarly, we obtain that$$[[x,e],J(x,y)]+[[x,y],J(x,e)]=0.$$
 Next, replacing $x$ by $x+e$ and $y$ by $y+e$ in the relation $[[x,y],J(x,y)]=0$ and summarizing the above conclusions, we see that $$[[x,y],J(e,e)]+[[x,e],J(e,y)]+[[e,y],J(x,e)]=0.$$ Hence, $[[exe,eye],J(e,e)]=0=[[e'xe',e'ye'],J(e,e)]$.\\ Also, from the relations $[[x,e],J(e,e)]=0$ and $e'J(e,e)e'=0$, we get $$J(e,e)exe'=[J(e,e),exe']e'=[[exe',e],J(e,e)]e'=0$$ and $$e'xeJ(e,e)=e'[e'xe,J(e,e)]=e'[[e'xe,e],J(e,e)]=0.$$ In a similar manner, we can show that $exe'J(e,e)=0=J(e,e)e'xe$.\\
 Therefore,  we conclude that $J_1$ is an extremal biderivation and $J_1(e,e)=J(e,e)$. Indeed,
\begin{eqnarray*}
[[x,y],J(e,e)]&=&[e[x,y]e+e[x,y]e'+e'[x,y]e+e'[x,y]e',J(e,e)]\\
&=&[e[x,y]e+e'[x,y]e',J(e,e)]\\
&=&[[exe,eye]+[e'xe',e'ye'],J(e,e)]\\
&=&0.
\end{eqnarray*}
It is easy to verify that $J_2=J-J_1$ is a Jordan biderivation.
\cqfd

 It is worthwhile mentioning that it was Herstein  who initiated the study of Jordan derivations on associative rings.  In  \cite{H}, he proved that for every   prime ring $A$ of characteristic different from 2  a Jordan derivation of $A$ is   a derivation of $A$. The following remark is the  ``biderivation" counterpart of his results  \cite[Lemmas 3.1.  and  3.2.]{H}.
 
\begin{remark}\label{rem-jbi}
    If $J:\mathcal{A}\times \mathcal{A}\longrightarrow \mathcal{A}$  is a Jordan biderivation, then the following assertions hold  for all $x,y,z, t\in \mathcal{A}$:
\begin{enumerate}
    \item  $J(xyx,z)=J(x,z)yx+xJ(y,z)x+xyJ(x,z)$.
 \item   $J(xyz+zyx,t)=J(x,t)yz+xJ(y,t)z+xyJ(z,t)+J(z,t)yx+zJ(y,t)x+zyJ(x,t)$.
\end{enumerate}
\end{remark}

The following lemma  is the key result for   decomposing a Jordan  biderivation as a sum of a biderivation and an antibiderivation.

\begin{lemma}\label{prop-jbider}
 Let $J:\mathcal{A}\times \mathcal{A}\longrightarrow \mathcal{A}$  be a Jordan  biderivation such that $J(e,e)=0$. Then, the following assertions hold
  for all $a,a'\in \mathcal{A}_{11},$ $ m,m'\in \mathcal{A}_{12},$ $ n,n'\in \mathcal{A}_{21}$ and $b,b'\in \mathcal{A}_{22}$:
 \begin{enumerate}
\item $J(a,a')=eJ(a,a')e $ and  $ J(b,b')=e'J(b,b')e'.$

\item $J(a,m)=aJ(e,m)+J(e,m)a$ and  $ J(m,a)=aJ(m,e)+J(m,e)a.$

\item $J(b,m)=bJ(e',m)+J(e',m)b$ and  $ J(m,b)=bJ(m,e')+J(m,e')b.$

\item $J(a,n)=aJ(e,n)+J(e,n)a$ and  $ J(n,a)=aJ(n,e)+J(n,e)a.$

\item $J(b,n)=bJ(e',n)+J(e',n)b$ and  $ J(n,b)=bJ(n,e')+J(n,e')b.$

\item $J(m,n)=eJ(m,n)e'+e'J(m,n)e+[J(e,n),m]=eJ(m,n)e'+e'J(m,n)e+[n,J(m,e)].$

\item $J(n,m)=eJ(n,m)e'+e'J(n,m)e+[J(n,e),m]=eJ(n,m)e'+e'J(n,m)e+[n,J(e,m)].$

\item $J(n,n')=eJ(n,n')e'+e'J(n,n')e+[n',J(n,e)]=eJ(n,n')e'+e'J(n,n')e+[n,J(e,n')].$

\item $J(m,m')=eJ(m,m')e'+e'J(m,m')e+[J(e,m'),m]=eJ(m,m')e'+e'J(m,m')e+[J(m,e),m'].$

\item $J(a,b)=J(b,a)=0.$
\end{enumerate}
\end{lemma}
\pr
From the assumption  $J(e,e)=0$ we get  $J(e,e')=J(e',e)=J(e',e')=0$. We have from $J(x^2,y)=xJ(x,y)+J(x,y)x$, then, taking $x=e$ and $y=a$, we get $eJ(e,a)e =0$. Similarly, from $J(x,y^2)=yJ(x,y)+J(x,y)y$, we get $eJ(a,e)e=0$. Now, according to  Remark    \ref{rem-jbi} and identity $eJ(e,a)e=eJ(a,e)e=0$, we have
\begin{eqnarray*}
J(a,a')=J(eae,ea'e)&=&eaJ(e,ea'e)+eJ(a,ea'e)e+J(e,ea'e)ae\\
&=&ea(ea'J(e,e)+eJ(e,a')e+J(e,e)a'e)\\
&+&e(ea'J(a,e)+eJ(a,a')e+J(a,e)a'e)e\\
&+&(ea'J(e,e)+eJ(e,a')e+J(e,e)a'e)ea\\
&=&eJ(a,a')e.
\end{eqnarray*}
Similarly, we can show that $J(b,b')=e'J(b,b')e'$ and that the relation $(10)$ is true. Now to prove the assertion (2), we use Remark    \ref{rem-jbi}. Namely,  we  obtain
\begin{eqnarray*}
J(a,m)&=&eaJ(e,eme'+e'me)+eJ(a,eme'+e'me)e+J(e,eme'+e'me)ae\\
&=&eaeJ(e,m)e'+e'J(e,m)eae\\
&=&aJ(e,m)+J(e,m)a.
\end{eqnarray*}
In a similar manner, we can prove that the  conditions $(3)$, $(4)$ and  $(5)$ hold. Next, we show that
$J(m,n)=eJ(m,n)e'+e'J(m,n)e+[J(e,n),m]=eJ(m,n)e'+e'J(m,n)e+[J(m,e),n]$, and one can prove analogously that the  conditions  $(7)$, $(8)$ and $(9)$ also hold. Indeed,
\begin{eqnarray*}
J(m,n)=J(eme'+e'me,n)&=&emJ(e',n)+eJ(m,n)e'+J(e,n)me'\\
&+&e'mJ(e,n)+e'J(m,n)e+J(e',n)me\\
&=&eJ(m,n)e'+e'J(m,n)e+[J(e,n),m].
\end{eqnarray*}
 On the other hand $J(m,n)=J(m,ene'+e'ne)=eJ(m,n)e'+e'J(m,n)e+[n,J(m,e)]$.
\cqfd

Consider the decomposition  of a bilinear mapping
$J:\mathcal{A}\times \mathcal{A}\longrightarrow \mathcal{A}$  given in (\ref{linear3}). When $J$ is a Jordan biderivation, we could continue this decomposition using all the assertions of  Lemma \ref{prop-jbider} so we get a new larger   decomposition.   We will show, in the following two lemmas, that one part $\Delta$ of this new decomposition is an antibiderivation and another part $D$ is a biderivation. So we get $$J=\Delta + D + eJ(m,n)e'+e'J(m,n)e+eJ(n,m)e'+e'J(n,m)e.$$
 Under the condition given in Theorem \ref{theorem-princ1}, we will show that  the   part $eJ(m,n)e'+e'J(m,n)e+eJ(n,m)e'+e'J(n,m)e$ is zero.

\begin{lemma}\label{lem-antibideri}
For a    Jordan biderivation $J:\mathcal{A}\times \mathcal{A}\longrightarrow \mathcal{A}$ such that $J(e,e)=0$, a mapping $\Delta:\mathcal{A}\times \mathcal{A}\longrightarrow \mathcal{A}$ is     an antibiderivation if it  satisfies the
   following conditions  for all $a,a'\in \mathcal{A}_{11},$ $m,m'\in \mathcal{A}_{12},$ $n,n'\in \mathcal{A}_{21}$ and $b,b'\in \mathcal{A}_{22}$:
   \begin{enumerate}
\item $\Delta(m,m')=e'J(m,m')e $ and $\Delta(n,n')=eJ(n,n')e'.$

\item $\Delta(a,m)=J(e,m)a$ and $\Delta(m,a)=J(m,e)a.$

\item $\Delta(b,m)=bJ(e,m)$ and $\Delta(m,b)=bJ(m,e).$

\item $\Delta(a,n)=aJ(e,n)$ and $\Delta(n,a)=aJ(n,e).$

\item $\Delta(b,n)=J(e,n)b$ and $\Delta(n,b)=J(n,e)b.$

\item $\Delta(a,a')=\Delta(b,b')=\Delta(b,b')=\Delta(a,b)=\Delta(b,a)=\Delta(m,n)= \Delta(n,m)=0.$
\end{enumerate}
\end{lemma}
\pr
Since $J(m,e)[A,A]=0$ and $J(a,m)=J(e,am)$, we conclude that $$\Delta(m,aa')-a'\Delta(m,a)-\Delta(m,a')a=J(m,e)aa'-a'J(m,e)a-J(m,e)a'a=0$$
and $$\Delta(a'm,a)-\Delta(m,a)a'-m\Delta(a,a')=J(m,a')a-J(m,e)aa'=J(m,e)[a',a]=0.$$
Using the fact  that $J(e,mb)=J(e,mb+bm)=J(e,m)b+bJ(e,m)=-J(b,m)$ and the condition $(3)$ in  Lemma \ref{prop-jbider}, we get
$$\Delta(a,mb)-\Delta(a,b)m-b\Delta(a,m)=J(e,mb)a-J(mb,e)a-b(J(e,m)a+J(m,e)a)=0.$$
 The other conditions for $\Delta$ to be an antibiderivation can be proved with a similar calculation.
\cqfd

In what follows we  will also use  Benkovi\v{c} and \v{S}irovnik's conditions \cite{BS}; that is   the algebra $\mathcal{A}$ satisfies the following two implications which will be refereed as ``the conditions \textbf{(*)}":
\begin{itemize}
    \item For  all $x \in  \mathcal{A}$, $exe \cdot e\mathcal{A} e' = \{0\} = e' \mathcal{A}e \cdot exe$ implies $ exe = 0$.
    \item For  all $x \in  \mathcal{A}$, $ e\mathcal{A} e' \cdot e' x e' = \{0\} = e' x e' \cdot e' \mathcal{A}e $ implies $e' x e' = 0$.
\end{itemize}

Some important examples of unital algebras with nontrivial idempotents having the conditions  \textbf{(*)}   are triangular algebras, matrix algebras and prime (and hence in particular simple) algebras with nontrivial idempotents.

\begin{lemma}\label{lem-bideri}
   Assume that  $\mathcal{A}$ satisfies the conditions  \textbf{(*)}.   Then, for a    Jordan biderivation $J:\mathcal{A}\times \mathcal{A}\longrightarrow \mathcal{A}$ such that $J(e,e)=0$,   a mapping $D:\mathcal{A}\times \mathcal{A}\longrightarrow \mathcal{A}$  is  a  biderivation  if it satisfies the    following conditions for all $a,a'\in \mathcal{A}_{11}, m,m'\in \mathcal{A}_{12}, n,n'\in \mathcal{A}_{21}$ and $b, b'\in \mathcal{A}_{22}$:
 \begin{enumerate}
\item $D(a,a')=J(a,a')$ and $D(b,b')=J(b,b').$

\item $D(a,m)=aJ(e,m)$ and $D(m,a)=aJ(m,e).$

\item $D(b,m)=J(e',m)b$ and $D(m,b)=J(m,e')b.$

\item $D(a,n)=J(e,n)a$ and $D(n,a)=J(n,e)a.$

\item $D(b,n)=bJ(e',n)$ and $D(n,b)=bJ(n,e').$

\item $D(m,m')=eJ(m,m')e'$ and $D(n,n')=e'J(n,n')e.$

\item $D(a,b)=D(b,a)=D(m,n)= D(n,m)=0.$
\end{enumerate}

\end{lemma}
\pr
From the assumption  $J(e,e)=0$ we get  $J(e,e')=J(e',e)=J(e',e')=0$ and as in the first part of the proof of Lemma \ref{prop-jbider},   $J(x^2,y)=xJ(x,y)+J(x,y)x$  and $J(x,y^2)=yJ(x,y)+J(x,y)y$, shows that  $eJ(e,a)e=eJ(a,e)e=0$. From    Remark    \ref{rem-jbi} and the identity $eJ(e,a)e=eJ(a,e)e=0$, we have
\begin{eqnarray*}
J(am,e)=J(am+ma,e)&=&J(eam+mae,e)\\
&=&J(e,e)am+eJ(a,e)m+aJ(m,e)+J(m,e)a\\
&+&mJ(a,e)e+maJ(e,e)e\\
&=&aJ(m,e)+J(m,e)a
\end{eqnarray*}
Using  condition  $(2)$ of Lemma \ref{prop-jbider} we get   $J(am,e)=J(m,a)$. Then, from $J(am,e)=aJ(m,e)+J(m,e)a$,  we get $eJ(m,e)e=0$. Indeed, we have $J(m,e)=J(em,e)=eJ(m,e)+J(m,e)e$. 
 Next, we claim that $J(a,a')m=[[a',a],J(m,e)]$. Since,
 \begin{eqnarray*}
a'aJ(m,e)+J(m,e)aa'&=&a'J(m,a)+J(m,a)a'\\
&=&a'J(am,e)+J(am,e)a'\\
&=&J(am+ma,a')\\
&=&J(a,a')m+aJ(m,a')+J(m,a')a\\
&=&J(a,a')m+aa'J(m,e)+J(m,e)a'a.
\end{eqnarray*}
Next, we claim that $J(a,a')m=[a',a]J(m,e)$ and $J(m,e)[a,a']=0$. Since,
 \begin{eqnarray*}
a'aJ(m,e)+J(m,e)a'a=J(m,a'a)&=&J(a'am,e)\\
&=&J(am,a')\\
&=&J(am+ma,a')\\
&=&aa'J(m,e)+J(a,a')m+J(m,e)a'a.
\end{eqnarray*}
 Now, using the conditions \textbf{(*)}, we get  $D(a_1a_2,a)=D(a_1,a)a_2+a_1D(a_2,a)$ for all $a_1, a_2, a \in \mathcal{A}_{11}$. Indeed,  for all $m\in \mathcal{A}_{12}$,  we have
 \begin{eqnarray*}
(D(a_1a_2,a)-a_1D(a_2,a)-D(a_1,a)a_2)m&=&(J(a_1a_2,a)-a_1J(a_2,a)-J(a_1,a)a_2)m\\
&=&[a,a_1a_2]J(m,e)-a_1[a,a_2]J(m,e)-[a,a_1]J(m,a_2)\\
&=&aa_1a_2J(m,e)-a_1a_2aJ(m,e)-a_1aa_2J(m,e)\\
&+&a_1a_2aJ(m,e)-aa_1J(m,a_2)+a_1aJ(m,a_2)\\
&=&aa_1a_2J(m,e)-a_1a_2aJ(m,e)-a_1aa_2J(m,e)\\
&+&a_1a_2aJ(m,e)-aa_1 a_2 J(m,e)+a_1aJ(m,e)a_2\\
&-&aa_1 J(m,e)a_2+a_1a a_2 J(m,e)\\
&=&[a_1,a]J(m,e)a_2\\
&=&J(a,a_1)ma_2\\
&=&0.
\end{eqnarray*}

Similarly we obtain, for all $n\in \mathcal{A}_{21}$,
$n(D(a_1a_2,a)-a_1D(a_2,a)-D(a_1,a)a_2)=0.$ So the first implication of the conditions \textbf{(*)}, gives $D(a_1a_2,a)=D(a_1,a)a_2+a_1D(a_2,a)$.\\
Moreover,
\begin{eqnarray*}
D(am,a')-aD(m,a')-D(a,a')m&=&a'J(am,e)-aa'J(m,e)-J(a,a')m\\
&=&a'aJ(m,e)-aa'J(m,e)-J(a',a)m\\
&=&0.
\end{eqnarray*}
Analogously,  we can prove the other relations for $D$ to be a biderivation.
\cqfd

We are now in a position to state and to prove the  second main result.

\begin{theorem}\label{theorem-princ1}
  Assume that  $\mathcal{A}$ satisfies the conditions \textbf{(*)}    and  that  the zero homomorphism is the  only  $(e\mathcal{A}e,e'\mathcal{A}e')$-module morphism $f:e\mathcal{A}e'\longrightarrow e\mathcal{A}e'$
 such that
$e[\mathcal{A},\mathcal{A}]e \cdot f(e\mathcal{A}e')=f(e\mathcal{A}e') \cdot e'[\mathcal{A},\mathcal{A}]e'=0$, then every Jordan  biderivation $J:\mathcal{A}\times \mathcal{A}\longrightarrow \mathcal{A}$   can be written as a sum of a biderivation, an antibiderivation and  an extremal biderivation.
\end{theorem}
\pr Let $J : \mathcal{A}\times \mathcal{A} \longrightarrow \mathcal{A}$ be a Jordan biderivation. Using  Lemma \ref{lemm-jbider2}, we get $J = J_1 + J_2$, where $J_1(x, y) = [x, [y, J(e, e)]]$ is an  extremal biderivation and $J_2$ is a Jordan biderivation with $J_2(e, e) = 0$.  Following the discussion given before Lemma \ref{lem-antibideri}, we get the result if we prove    that  $$eJ(m,n)e'=e'J(m,n)e=eJ(n,m)e'=e'J(n,m)e=0.$$
Now fix $n\in \mathcal{A}_{21}$. Then, the map $f:\mathcal{A}_{12}\longrightarrow \mathcal{A}_{12}$ defined by $f(m)=eJ(m,n)e'$  (for all $m\in \mathcal{A}_{12}$) is a module homomorphism. Moreover, according to   Remark    \ref{rem-jbi},  we get $[\mathcal{A}_{11},\mathcal{A}_{11}]f(m)=f(m)[\mathcal{A}_{22},\mathcal{A}_{22}]=0$. Indeed,   for all $a,a'\in \mathcal{A}_{11}$ and $b,b'\in \mathcal{A}_{22}$,  one can check easily that  $aa'f(m)=a'af(m) $
  and  $f(m)bb' =f(m)b'b.$
  Hence,   by hypothesis,  $f=0$. Similarly, we can obtain the other relations.\cqfd

Let $Id([\mathcal{A},\mathcal{A}])$ denotes the ideal generated by all commutators $[x, y]$ ($x, y\in \mathcal{A}$) of an algebra $\mathcal{A}$.
The following corollary is  an immediate consequence of Theorem \ref{theorem-princ1}.

 \begin{corollary}\label{coroll-princ1}   If    $\mathcal{A}$ satisfies the conditions \textbf{(*)}  and   either $Id([e\mathcal{A}e,e\mathcal{A}e])=e\mathcal{A}e$ or $Id([e'\mathcal{A}e',e'\mathcal{A}e'])=e'\mathcal{A}e'$, then every Jordan  biderivation  $J:\mathcal{A}\times \mathcal{A}\longrightarrow \mathcal{A}$  can be written as a sum of a biderivation and an antibiderivation.
 \end{corollary}
\pr
We have $J(m,e)[\mathcal{A}_{11},\mathcal{A}_{11}]=0$, then by hypothesis, we get $J(m,e)\mathcal{A}_{11}=0$. Since $\mathcal{A}_{11}$ has a unity element, it follows that $J(m,e)=0$. So $J(a,a')m=[a',a]J(m,e)=0$ for all $m\in \mathcal{A}_{12}$. Similarly we have $nJ(a,a')=0$ for all $n\in \mathcal{A}_{21}$. Thus, using the conditions \textbf{(*)},  $J(e,e)=0$. Therefore, by    Theorem \ref{theorem-princ1}, we get the result. \cqfd


If $\mathcal{A}$ admits a nontrivial idempotent $e$ such that $e\mathcal{A}  e'\mathcal{A}e = \{0\} = e'\mathcal{A}  e\mathcal{A}e'$  and the bimodule $e\mathcal{A}f$ is faithful as both a left $e\mathcal{A}e$-module and     a right $e'\mathcal{A}e'$-module, then $\mathcal{A}$ satisfies the conditions  \textbf{(*)}. Then, we have  the following corollary.

\begin{corollary}
Assume that  $e\mathcal{A}e'  \mathcal{A}e = \{0\} = e'\mathcal{A}  e\mathcal{A}e'$  and  either   $Id([e\mathcal{A}e,e\mathcal{A}e])=e\mathcal{A}e$ or $Id([e'\mathcal{A}e',e'\mathcal{A}e'])=e'\mathcal{A}e'$. If the bimodule $e\mathcal{A}f$ is faithful as both a left $e\mathcal{A}e$-module and     a right $e'\mathcal{A}e'$-module,  then every Jordan  biderivation  $J:\mathcal{A}\times \mathcal{A}\longrightarrow \mathcal{A}$   can be written as a sum of a biderivation and an antibiderivation.
\end{corollary}

Recall that   an algebra  $\mathcal{A}$ is isomorphic to a  triangular matrix algebra     if  there exists a non trivial idempotent $ e \in \mathcal{A}$ such that $ (1-e)\mathcal{A}e=0 $ (see, for instance, \cite[Theorem 5.1.4]{GJS}). Thus,  triangular algebras  are examples  of algebras that satisfies the  conditions of Theorem \ref{theorem-princ1}. Namely, we get the following result which generalizes \cite[Theorem 2.10]{AHA}.

 \begin{corollary}\label{theorem-princ2}
 Every Jordan biderivation on a triangular algebra is a biderivation.
 \end{corollary}

 At this stage we remark that, triangular algebras are special examples of trivial extension algebras on which the Jordan generalized and Lie generalized derivations are recently investigated  in \cite{BBEEF, BEFKM}. We conclude this section with the following question, to the best of our knowledge, has not been studied/answered  yet:
 \begin{question}
 Under what conditions every Jordan biderivation on a trivial extension algebra is a biderivation?
 \end{question}

\noindent\textbf{Acknowledgement.} The authors would like to thank the referee for the careful reading of the paper.

\begin{center}\small{1.   Department of Pure Mathematics, Ferdowsi University of Mashhad, P.O. Box 1159, Mashhad 91775, IRAN.}\end{center}\vspace{-0,4cm}

\noindent\hspace{4cm}\small{mohamadali$_-$bahmani@yahoo.com}\smallskip

\begin{center}
\small{2.  Centre de Recherche de Math\'ematiques et Applications de Rabat (CeReMAR), Faculty of Sciences, Mohammed V University in Rabat,  Morocco.}
\end{center}\vspace{-0,4cm}

\noindent\hspace{4cm}\small{driss.bennis@um5.ac.ma; driss$\_$bennis@hotmail.com}\smallskip

\begin{center}\small{3.   Department of Pure Mathematics, Centre of Excellence in Analysis on Algebraic Structures (CEAAS), Ferdowsi University of Mashhad, P.O. Box 1159, Mashhad 91775, IRAN.}\end{center}\vspace{-0,4cm}

\noindent\hspace{4cm}\small{vishki@um.ac.ir}\smallskip

\begin{center}
\small{4.  Superior School of Technology, Ibn Tofail University, Kenitra,  Morocco.}
\end{center}\vspace{-0,4cm}

\noindent\hspace{4cm}\small{fahid.brahim@yahoo.fr} 
\end{document}